\documentclass[12pt,twoside]{article}
\usepackage{latexsym,amsfonts,amssymb,mathrsfs,amsthm}
\usepackage[leqno]{amsmath}
\numberwithin{equation}{section}
\usepackage{color}
\usepackage{colordvi,multicol}
\usepackage{fancyhdr}
\DeclareMathOperator{\var}{var}
\newtheorem{thm}{\hskip\parindent\bf{Theorem}}[section]
\newtheorem{cor}{\hskip\parindent\bf{Corollary}}[section]
\newtheorem{lem}{\hskip\parindent\bf{Lemma}}[section]

\newtheorem{rem}{\hskip\parindent\bf{Remark}}[section]
\newtheorem{exa}{\hskip\parindent\bf{Example}}[section]

\makeatletter
\def\@makefnmark{}
\makeatother

\begin{document}
\title{\bf Exponential bounds of  ruin  probabilities   for non-homogeneous  risk models}
\author{\bf \small{BY} \\   \bf \small{QIANQIAN ZHOU, ALEXANDER SAKHANENKO}\\ \bf \scriptsize{AND} \bf \small{JUNYI GUO$^{*}$}
 }
\date{}
\maketitle
 { Abstract.}\quad
 Lundberg-type inequalities for ruin probabilities  of   non-homogeneous risk models  are  presented  in this paper. By employing  martingale method,   the  upper bounds of ruin probabilities   are obtained  for the  general  risk  models  under  weak    assumptions.  In addition, several  risk  models,  including
 the newly defined united risk model and quasi-periodic risk model with  interest rate,   are studied.



\footnote{The research of Qianqian Zhou and Junyi Guo was supported by the National Natural Science
Foundation of China (Grant No. 11931018) and Tianjin Natural Science Foundation.

The research of Alexander Sakhanenko  was also supported by the program of fundamental scientific researches of the SB RAS № I.1.3., project № 0314-2019-0008.

$^*$ Corresponding  author.}
\thispagestyle{empty}
\section{Introduction}
Research  on the ultimate ruin  probability $\psi(u)$, i.e., the probability that the reserve ever drops below zero:
\begin{gather*}                                                       
\psi(u):={\bf P}[\inf_{t\geq 0}R(t)<0| R(0)=u],
\end{gather*}
where   $R(t)$ is a risk reserve process with  initial reserve $ R(0)=u>0,$
 is attracting  increasing  attention  since classic  works  of  Lundberg \cite{LF03} and Cram$\acute{e}$r \cite{C30,C55}.
  After that, a  substantial  amount of works  have been devoted to finding   the ruin probabilities  and    the upper  bounds  of ruin probabilities, see, e.g.,      Dickson \cite{DC05}, Gerber \cite{GH79}, and Rolski et al. \cite{RSST98}.
However, the  most outstanding result about the behavior of $\psi(u)$
 is the Lundberg  inequality,  which states that under appropriate assumptions (see  \cite{AA10}  and \cite {G91}  for  more details)
\begin{gather}                                                                                \label{a2}
\psi(u)\leq e^{-Lu} \text{ for  all } u\geq 0.
\end{gather}
The largest number $L$ is unique and called the adjustment coefficient or Lundberg exponent.
Further,    the following   Lundberg-type  inequality is  also  studied
\begin{gather}                                                        \label{a3}
\text{ for all } u\ge0,\  \psi(u)\leq Ce^{-Lu}\text{ with }  C<\infty\text{ and }  L>0.
\end{gather}

 In  the  classical  works  of   homogeneous  risk  models, both inter-occurrence times and claim sizes are assumed to be  i.i.d. random variables.  However, in reality, both factors are influenced by the economic environment. For instance, inflation and interest rate can affect the evolution of the reserves of the company. Thus  the assumptions of homogeneity of
inter-occurrence times  and claim amounts can be too restrictive  for practical use.  It is obvious  that  the non-homogeneous models  better reflect  the real insurance activities compared  to  the homogeneous risk models.

Recently, the attention on non-homogeneous risk models has considerably increased.  In the non-homogeneous risk models the inter-occurrence times or (and) the claim sizes are independent but not necessarily identically distributed. In  Bla$\check{z}$evi$\check{c}$ius,  et  al. \cite{BBE10}, Casta$\tilde{n}$er, et al. \cite{CCGLM13}, Lef$\grave{e}$vre  and Picard \cite{LP06},     and  Vernic \cite{VR15}, they investigated the non-homogeneous risk models with independent, identically distributed inter-occurrence times but not necessarily identically distributed claims.   For more details about  the non-homogeneous risk models with  claim sizes, which are independent and identically distributed random variables, and independent but not identically distributed inter-occurrence times one can refer to Bernackait$\dot{e}$ and $\breve{S}$iaulys  \cite{BJ17}, Ignator  and Kaishev \cite{IK00} and  Tuncel and  Tank \cite{TT14}.

But we found only a few works in which estimates of the {types \eqref{a2} or \eqref{a3}
were obtained.
In  Andrulyt$\dot{e}$, et al. \cite{AE15},  Kievinait$\dot{e}$  and $\check{S}$iaulys \cite{KD18}, and Kizinevi$\check{c}$  and $\check{S}$iaulys \cite{EJ18},   the non-homogeneous renewal risk model, where claim sizes and  inter-occurrence times are both independent but not necessarily identically distributed, is considered and   the Lundberg type inequality of this model is investigated.

This paper considers the  non-homogeneous compound Poisson risk  model  and   studies its   ruin  probability  $\psi(u).$
Our  purpose  is   to   obtain   the  exponential  upper  bound  of  ultimate  ruin  probability   under  general  assumptions.

The  rest  of  the  paper is  organized  as  follows.  In   Section \ref{section2}, we  present  our  main  results. First, we study the non-homogeneous compound Poisson risk model without interest rate. By employing martingale method, the Lundberg inequality of this model is obtained. Second,  a  non-homogeneous risk model with interest rate is  studied   and   the  Lundberg  inequality  of  its    ruin probability     is   also  obtained.
In Section \ref{section3}, the applications of our main results with  a special attention  to new possibilities of the behavior of probability $\psi(u)$ are considered.   First,  the periodic and quasi-periodic risk models with interest rate, which are also automatically non-homogeneous, are investigated. Second,  a new non-homogeneous risk model, which is  called united risk model, is defined. And  the Lundberg inequality of its probability of ruin is also obtained. In the last section, we compare our results with the previous results in Andrulyt$\dot{e}$, et al. \cite{AE15},  Kievinait$\dot{e}$  and $\check{S}$iaulys \cite{KD18}, and Kizinevi$\check{c}$  and $\check{S}$iaulys \cite{EJ18}.


\section{Main results}\label{section2}

\subsection{Non-homogeneous risk model without interest rate}\label{sec2.1}

We begin to consider   a  general  risk model in which the distribution of each
inter-arrival time $\theta_n:= T_n-T_{n-1}$ may depend on time $T_{n-1}$.

{\bf Model  A. }Let $N(t)$ be a non-homogeneous  Poisson process with non-decreasing
accumulate  intensity function $\Lambda(t):={\bf E}N(t)<\infty$ where $N(0)=\Lambda(0)=0$.
Denote by $Z(t)\geq 0$  the claim size in case when it  arrived at time $t>0$ and   for each  j=0,1,2,\dots,  let
$$T_j:=\min\{t\ge0: N(t)\geq j\}.
$$
 Then   it is natural to introduce risk reserve process $R(t)$ in the following way
\begin{gather*}                                                                         
R(t)= u+p(t)-\sum_{j=1}^{N(t)}Z(T_j), \  \forall t\geq 0,
\end{gather*}
where $u=R(0)>0$  is  the  initial  reserve,  $p(t)\geq 0$  is the aggregate premium income over $[0, t],$   which is assumed to be non-decreasing   and   $\{Z(T_j)\}_{j\ge 1}$  are the   claim    sizes which are  supposed  to  be  independent  but   not  identically  distributed.

It  is  easy  to  see  that  model  A  is   more  general  than  the  classical  compound  Poisson  risk  model which  supposes   that  $N(t)$  is  the homogeneous   Poisson  process, $\{Z(T_j)\}_{j\ge 1}$  are  i.i.d. random variables  and  the  premium  rate is a  constant.

For  mathematical  purpose, introduce now the claim surplus process
\begin{gather}                                                                         \label{d2}
S(t) 
=\sum_{j=1}^{N(t)}Z(T_j)-p(t),
\end{gather}
thus $R(t)=u-S(t).$
Then   the  ultimate  ruin  probability  can  be  expressed  as
\begin{gather*}
\psi(u)={\bf P}(\sup_{t\ge 0}S(t)>u).
\end{gather*}

We suppose that random processes $(N(\cdot), Z(\cdot))$  are independent. From Kingman \cite{KJF93}   and  Paulsen \cite{PJ97} we can see that  for any real number $h\ge0$ and any $t>0$  if $Q(h, t)={\bf E}e^{hZ(t)}-1$  exists and is continuous  in $t\in(0, \infty),$  then for the  claim surplus process $S(t)$  defined in (\ref{d2})
\begin{gather}                          \label{d3}
{\bf E}e^{hS(t)}=\exp\Big\{ \int_{0}^{t} Q(h, x)d\Lambda(x)-hp(t)     \Big\},\  \forall h\ge 0,\   t>0,
\end{gather}
where $\Lambda(t):={\bf E}N(t)$  is the accumulated intensity function of $N(t)$ and $p(t)$ is the  aggregate premium income over $[0, t].$
Note that the formula in (\ref{d3}) is valid for a larger class of non-homogeneous  processes $S(t)$  with independent increments.

Before studying model A,  we  first  state    the  following  basic  theorem
\begin{thm}\label{thm0}
Let   $W(t)$  be     a separable process  with independent increments, then for any real numbers $x$  and $T, h\ge 0$
\begin{gather}\label{01}
  {\bf P}(\sup_{0\le t\le T}W(t)>x)\le  e^{-hx}\sup_{0\le t\le T}{\bf E}e^{hW(t)}.
\end{gather}
In addition, for  all $x$  and $h\ge 0$
\begin{gather}\label{02}
 {\bf P}(\sup_{t\ge 0}W(t)>x)\le  e^{-hx}\sup_{t\ge 0}{\bf E}e^{hW(t)}.
\end{gather}
\end{thm}

Before proving Theorem \ref{thm0}  we first introduce a key lemma which comes from paper \cite{ZSG19}.
\begin{lem} 
\label{L-1}
If  random  variables $Y_1, Y_2, \ldots$  are mutually independent, then for any $n\ge 1$ and any real number $x$ and $h\ge 0$
\begin{gather*}
 {\bf P}[\sup_{1\le k\le n}W_k>x]\leq e^{-hx}\sup_{1\le k\le n}{\bf E}e^{h W_k},
\end{gather*}
where  $W_k:=Y_1+Y_2+\ldots+Y_k.$  In addition, for any real number $x$ and $h\ge 0$
\begin{gather*}
{\bf P}[\sup_{k\ge 1}W_k>x]\leq e^{-hx}\sup_{k\ge 1}{\bf E}e^{h W_k},
\end{gather*}
\end{lem}

\noindent{\bf  Proof  of  Theorem \ref{thm0}.}
Since $W(t)$  is a separable  random process, then  there  exists  a sequence  $ t_1, t_2, \ldots\in [0, T]$  such that
$$\max_{1\leq k\leq n}W(t_k)\uparrow \sup_{k\geq 1} W(t_k)=\sup_{0\le t\le T}W(t).$$
Thus
$${\bf P}(\sup_{0\le t\le T}W(t)>x)=\lim_{n\rightarrow \infty}{\bf P}(\max_{1\leq k\leq n}W(t_k)>x).$$

Let $\{t_1, t_2, \ldots, t_n\}=\{ t^{*}_1<t^{*}_2<\ldots< t^{*}_n    \},$  then
$$\max_{1\leq k\leq n}W(t_k)=\max_{1\leq k\leq n}\sum_{j=1}^{k}[W(t^{*}_j)-W(t^{*}_{j-1})],$$
where $t^{*}_0=0,\ W(t^{*}_0)=0.$
Since $\{W(t^{*}_j)-W(t^{*}_{j-1})\}_{j\geq 1}$  are independent,  then  from   Lemma \ref{L-1} we have
\begin{gather}\label{m3}
{\bf P}(\max_{1\leq k\leq n}W(t_k)>x)\leq e^{-hx}\max_{1\leq k\leq n}{\bf E}e^{hW(t_k)}.
\end{gather}
Let $n\rightarrow \infty $  on both  sides of (\ref{m3}), we have
$${\bf P}(\sup_{0\le t\le T}W(t)>x)
\leq e^{-hx}\sup_{0\le t\le T}{\bf E}e^{hW(t)}.$$

Let  $T\rightarrow \infty$  in (\ref{01})  then (\ref{02}) is  obtained.

\hfill\fbox

Note  that    Theorem \ref{thm0} allows  us to  obtain  the  following  results  of  the ruin  probabilities  of model A.
\begin{thm}\label{thm-1}
Suppose
the  claim  surplus  process   $S(t)$  defined  in  (\ref{d2}) is a separable process  with independent increments.
Then for any $u>0$ and $T, h\ge0$
\begin{gather*}                                                                        
  \psi(u, T)= {\bf P}(\sup_{0\le t\le T}S(t)>u)\leq e^{-hu}\sup_{0\le t\le T}{\bf E}e^{hS(t)}.
\end{gather*}
In addition, for any $u>0$ and $h\ge0$
\begin{gather*}
  \psi(u)= {\bf P}(\sup_{t\ge 0}S(t)>u)\leq e^{-hu}\sup_{t\ge 0}{\bf E}e^{hS(t)}.
\end{gather*}
\end{thm}

 Theorem \ref{thm-1}  and (\ref{d3}) immediately imply  the following Lundberg inequality of ruin probability of  model A.
\begin{cor}\label{c1}
Under the conditions of Theorem \ref{thm-1},
 if  $S(t)$   satisfies  (\ref{d3}) then   the Lundberg   inequality (\ref{a2})  holds with $L=\underline{L},$
where
\begin{gather*}                                                                        
\underline{L}=\sup\Big\{h\geq0:       \sup_{t\ge 0} \big[\int_{0}^{t}Q(h, x)d \Lambda(x)-hp(t)    \big]\leq 0    \Big\}.
\end{gather*}

\end{cor}

Homogeneous  compound  Poisson  risk  model  is a  special  case  of  model A, if
\begin{gather}                                                                         \label{d9}
{\bf E}e^{hS(t)}=\exp\Big \{t\big[\lambda({\bf E}e^{hZ}-1)-hp\big]     \Big \},\  \forall t>0,\  \forall h\ge 0.
\end{gather}
Then  from   Corollary \ref{c1}  the  classical   Lundberg inequality  (\ref{a2})  takes  place  with
\begin{gather}\label{d10}
L=\sup \Big\{h\ge 0 :\lambda({\bf E}e^{hZ}-1)\leq hp   \Big\}.
\end{gather}
\begin{rem}\label{r1}
From  (\ref{d10})  it  can  be   seen  that the  Lundberg inequality  in  our  method   is  better  than the famous Lundberg inequality in its  classical  form,  which  can  be  found in papers \cite{AA10,G91,RSST98},
 because
  we do not exclude
the cases when
${\bf E}e^{LS(t)}<1$ and/or when the expectation ${\bf E}S(t)=-\infty.$
\end{rem}

\subsection{ Non-homogeneous risk model with  interest rate}
As we  all  know, interest rate is one of the important factors that affect the evolution  of  the  reserves  of  the  company. And  underline that all risk models under rates of interest are automatically
non-homogeneous. This fact gives us additional motivations in our research of non-homogeneous risk
models. Thus   in  this  subsection  we  consider   the  non-homogeneous  compound Poisson risk model  with  interest  rate.

Here we    consider the  following  model

$\textbf{Model B: }$
Let  $N(t)$, $Z(t)$ and $p(t)$ be  as  in  model A, suppose that for any  $ t_2>t_1\ge0$ risk reserve process $R(t)$ satisfies  the following property
\begin{gather}\label{d11}
R(t_2)=(1+\alpha(t_1, t_2))R(t_1)+(1+\beta(t_1, t_2))[p(t_2)-p(t_1)]-[q(t_2)-q(t_1)],
\end{gather}
where  $(1+\alpha(t_1, t_2))R(t_1)$ is the accumulated value of $R(t_1)$ from $t_1$ to $t_2$ under interest rate  $\alpha(t_1, t_2)$
and $(1+\beta(t_1, t_2))[p(t_2)-p(t_1)]$ is the accumulated value of premiums collected from $t_1$ to $t_2$ under interest rate  $\beta(t_1, t_2)$.

Here  $q(t)$ denotes the total  claims over  time interval $[0, t]$. Of course, in the first order we will consider the case when $q(t)=\sum_{j=1}^{N(t)}Z(T_j)$.
We  suppose that
\begin{gather}\label{d12}
\beta(t_1, t_2)\geq 0
 \text{ and }
 \alpha(t_1, t_2)\geq e^{r(t_2)-r(t_1)}, \  \forall t_2>t_1\ge0,
\end{gather}
where $r(t)$  is non-decreasing and  $r(0)=0.$
So  we have   from  (\ref{d11})  and (\ref{d12})   that
\begin{gather*} 
R(t_2)\geq e^{r(t_2)-r(t_1)}R(t_1)-[S(t_2)-S(t_1)],
\end{gather*}
with   $S(t)=q(t)-p(t).$  Hence for  all   $t_2>t_1\geq0$
\begin{gather}\label{d13}
e^{-r(t_2)}R(t_2)\geq e^{-r(t_1)}R(t_1)-e^{-r(t_2)}[S(t_2)-S(t_1)].
\end{gather}

From  the inequality  (\ref{d13})  for  any  $0=t_0< t_1<t_2<\cdots <t_m = t$   we have
\begin{eqnarray}                                                                                              \label{d14}
e^{-r(t)}R(t)-u&=&\sum_{k=1}^{m}[e^{-r(t_k)}R(t_k)- e^{-r(t_{k-1})}R(t_{k-1})]\nonumber\\
&\geq&
-\sum_{k=1}^{m}e^{-r(t_k)}[S(t_k)-S(t_{k-1})].
\end{eqnarray}
If  $S(t)$  is a function of bounded   variation, then standard arguments allow us, as
$m\rightarrow \infty$, to change the sum  in the right hand side of  (\ref{d14})
by the integral
\begin{gather}                                                                                              \label{d15}
e^{-r(t)}R(t)\geq u-\int_{0}^{t}e^{-r(x)}d S(x)=u-Y(t),\  \forall  t>0.
\end{gather}
In addition, from (\ref{d15})we have
\begin{gather}                                                                                            \label{d16}
\psi(u)={\bf P}(\inf_{t\ge 0}R(t)<0)={\bf P}(\inf_{t\ge 0}e^{-r(t)}R(t)<0)\leq {\bf P}(\sup_{t\ge 0}Y(t)>u).
\end{gather}

It  is easy  to  see  that if the integral  in  (\ref{d15})  has sense  for   some
process   $S(t)$  with independent  increments and  non-random  $r(t),$
then $Y(t)$  is  also   a process    with independent  increments.
This fact and (\ref{d15}) imply the following assertion.

\begin{thm}\label{thm2}
Suppose that  (\ref{d15})   holds   where  $Y(t)$ is  a separable  process with independent  increments.
Then all  the  assertions  of    Theorem \ref{thm-1} take place   with   $Y(t)$  instead  of  $S(t).$
\end{thm}
\begin{cor}\label{c2}
Under  the  assumptions of  Theorem \ref{thm2}, also assume  that
\begin{gather}\label{thm2-1}
{\bf E}e^{hY(t)}=e^{a(h,t)}  \text{ with }
a(h,t):=-h\int_{0}^{t}e^{-r(x)}d p(x) +\int_{0}^{t}Q(he^{-r(x)}, x)d \Lambda(x).
\end{gather}
Then   the Lundberg   inequality (\ref{a2})  holds   with $L=L^*,$
where
$$L^{*}=\sup\big\{h\ge0:\sup_{t\ge 0}a(h,t)\le0\big\}.$$
\end{cor}

This assertion is an analogue of Corollary \ref{c1} for models
without  rates of interest.

\begin{rem}   \label{R-2}
Suppose   the  representation  (\ref{d3})  takes  place with a non-decreasing function $\Lambda(\cdot)$ and  natural continuity assumptions on function $Q(\cdot,\cdot)$ in the domain where it is finite.
It  is not  difficult to   understand that if the
 integral  in  (\ref{d15})  has sense  for   some
process   $S(t)$  with independent  increments and  non-random  function $r(t)$ of bounded variation   then  $Y(t)$  is  also   a process    with independent  increments such that
\begin{gather*}
{\bf E}e^{hY(t)}=e^{a(h,t)} \text{ with }
a(h,t):=-h\int_{0}^{t}e^{-r(x)}d p(x) +\int_{0}^{t}Q(he^{-r(x)}, x)d \Lambda(x).
\end{gather*}
We  can refer to Paulsen \cite{PJ97}  for more details. Thus  the  assumption  (\ref{thm2-1})  in  Corollary  \ref{c2}  is  reasonable
\end{rem}

\begin{rem}
Note  that  Theorems  \ref{thm-1}  and  \ref{thm2}  immediately  come  from  Theorem \ref{thm0}.
\end{rem}

\section{Applications of main results}\label{section3}
In this section, we  consider several applications of our general results with a special attention to
 new possibilities for the behavior of probability $\psi(u)$ in non-homogeneous cases. First, we investigate  periodic and quasi-periodic risk models which are also automatically non-homogeneous. Second, we put forward a new non-homogeneous risk model and study its probability of ruin.

\subsection{ Periodic and quasi-periodic risk models}
In  this  subsection,  we   study  the risk models which  take place in  the
periodic environment. For  the  classical  periodic  case,  Asmussen and  Rolski \cite{AR94}  has shown that the adjustment coefficient in the periodical
risk model is the same as for the standard time-homogeneous
compound Poisson risk process obtained by averaging the
parameters over one period. Here we study the general  periodic risk models with  interest rates.

In the  following,  under   the  conditions  of  Theorem  \ref{thm2}, we  also   assume  that
  $\Lambda(t)$  has density  $\Lambda^{'}(t)$ and
the aggregate premium income  $p(t)$   has density $p^{'}(t)$. Then (\ref{thm2-1})  can be  rewritten  as
\begin{gather*}  
 {\bf E}e^{hY(t)}=\exp\big\{ -\int_{0}^{t}he^{-r(x)}p^{'}(x)dx+\int_{0}^{t}({\bf E}e^{he^{-r(x)}Z(x)}-1)\Lambda^{'}(x)dx \big\},\    \forall  t>0.
\end{gather*}

\begin{cor}  \label{C-3}
Suppose there exist real  numbers  $l, t_0>0$  and   $\widetilde{L}>0$   such that for any $t\ge t_0$
\begin{gather}  \label{p1}
  {\bf E}e^{\widetilde{L}(Y(t+l)-Y(t))}\le 1.
\end{gather}
Then under the conditions of Theorem \ref{thm2}, the assertions in Theorem \ref{thm2} hold for each $h\in[0, \widetilde{L}]$ with
\begin{gather}  \label{p1+}
    \sup_{t\ge  0}{\bf E}e^{hY(t)}=\sup_{0\le t<  t_0+l}{\bf E}e^{hY(t)}.
\end{gather}
In addition, for any $u>0$
$$     \psi(u)\le \inf_{h\in [0,   \widetilde{L}]}\big\{ e^{-hu} \sup_{0\le t<  t_0+l}{\bf E}e^{hY(t)}   \big\}\le C_1 e^{-\widetilde{L}u},$$
where $C_1:=\sup_{0\le t<  t_0+l}{\bf E}e^{\widetilde{L}Y(t)}. $
\end{cor}
{\bf   Proof.}
For  any  $t\ge t_0+l,$   random  variables   $\Delta_{n,l}=Y(t)-Y(t-l)$   and  $Y(t-l)$   are independent.
Hence for  each  $h\in [0, \widetilde{L}],$
\begin{gather}   \label{m4}
{\bf E}e^{h\Delta_{n,l}}\le ({\bf E}e^{\widetilde{L}\Delta_{n,l}})^{h/\widetilde{L}}\le 1 \text{ and } {\bf E}e^{hY(t)}={\bf E}e^{h\Delta_{n,l}}{\bf E}e^{hY(t-l)}\le {\bf E}e^{hY(t-l)}.
\end{gather}
Using induction with respect to $t$ it is not difficult to obtain from \eqref{m4} that for any $ t\ge t_0+l$ and $h\in [0, \widetilde{L}]$
\begin{gather*}
 {\bf E}e^{hY(t)}\le \sup_{ t<t_0+l}{\bf E}e^{hY(t)}.
\end{gather*}
Hence
$$  \sup_{0\le t<t_0+l}{\bf E}e^{hY(t)}\le\sup_{t\ge 0}{\bf E}e^{hY(t)}\le \sup_{0\le t< t_0+l}{\bf E}e^{hS(t)}, \   \forall  h\in [0, \widetilde{L}],$$
i.e.,
\begin{gather*}
\sup_{t\ge 0}{\bf E}e^{hS(t)}=\sup_{0\le t< t_0+l}{\bf E}e^{hS(t)},\  \forall  h\in [0, \widetilde{L}].
\end{gather*}

The   rest of  assertions  come  from Theorem \ref{thm2} and Corollary \ref{c2}.

\hfill\fbox

\begin{thm}\label{thm3}
Suppose there exists a real  number $l>0$ such that  for each $h\in [0, L(Y(l))]$  and any $t\ge 0$
\begin{gather}                                           \label{p18}
   e^{r(t+l)} \big[{\bf E}e^{he^{-r(t+l)}Z(t+l)}-1\big]\leq e^{r(t)}\big[{\bf E}e^{he^{-r(t)}Z(t)}-1\big],
\end{gather}
where
\begin{gather}\label{p18+}
L(Y(l)):=\sup\big\{ h\ge 0:    {\bf E}e^{hY(l)}\le 1      \big\}.
\end{gather}
Assume also that for any $t\ge 0$
\begin{gather}                                           \label{p19}
 e^{-r(t+l)}\Lambda^{'}(t+l)\leq e^{-r(t)}\Lambda^{'}(t) \text{ and }
e^{-r(t+l)}p^{'}(t+l)\ge e^{-r(t)}p^{'}(t).
\end{gather}
Then under the conditions of Theorem \ref{thm2}, all the assertions in Theorem \ref{thm2} hold
   for each $h\in[0, L(Y(l))]$  with
\begin{gather} \label{p19+}
\sup_{t\ge  0}{\bf E}e^{hY(t)}= \sup_{0\le t\le l}{\bf E}e^{hY(t)}.
\end{gather}
In  addition,  for any $u>0$
\begin{gather*}                                          
 \psi(u)\leq  \inf_{h\in [0,  L(Y(l))]}\{e^{-hu}\sup_{0\le  t\leq l}{\bf E}e^{hY(t)}\}\le C_2 e^{-L(Y(l))u},
\end{gather*}
where $C_2:=\sup_{0\le t\le l}{\bf E}e^{L(Y(l))Y(t)}.$
\end{thm}
{\bf   Proof.}
For any    $t>l,$  random  variables  $\Delta_{n,l}:=Y(t)-Y(t-l)$ and $Y(t-l)$
are independent, then
\begin{eqnarray*}
{\bf E}e^{hY(t)}={\bf E}e^{hY(t-l)}{\bf E}e^{h\Delta_{n,l}}.
\end{eqnarray*}

Denote
$$F(t):=\int_{t-l}^{t}G(x)dx,\text{where }G(x):=-he^{-r(x)}p^{'}(x)+({\bf E}e^{he^{-r(x)}Z(x)}-1)\Lambda^{'}(x).$$
Then
$${\bf E}e^{h\Delta_{n,l}}=e^{F(t)}.$$

 Under the conditions of (\ref{p18})  and (\ref{p19})  we have
\begin{eqnarray*}
F^{'}(t)&=&G(t)-G(t-l)\nonumber\\
&=&  h(e^{-r(t-l)}p^{'}(t-l)-e^{-r(t)}p^{'}(t))+({\bf E}e^{he^{-r(t)}Z(t)}-1)\Lambda^{'}(t)
-\\&&({\bf E}e^{he^{-r(t-l)}Z(t-l)}-1)\Lambda^{'}(t-l)\nonumber\\
&=& h(e^{-r(t-l)}p^{'}(t-l)-e^{-r(t)}p^{'}(t))+e^{-r(t)}\Lambda^{'}(t)\cdot e^{r(t)}\big[ {\bf E}e^{he^{-r(t)}Z(t)}-1   \big]
-\\&&e^{-r(t-l)}\Lambda^{'}(t-l)\cdot e^{r(t-l)}\big[ {\bf E}e^{he^{-r(t-l)}Z(t-l)}-1   \big]\nonumber\\
&\le&  0.
\end{eqnarray*}
Thus  $F(t)$  is  non-increasing  such  that
$${\bf E}e^{h\Delta_{n,l}}=e^{F(t)}\le e^{F(l)}={\bf E}e^{hY(l)}\le 1,\  \forall  h\in [0, L(Y(l))]. $$
Therefore, we obtain that
\begin{gather}   \label{m5}
{\bf E}e^{hS(t)}\leq {\bf E}e^{hS(t-l)},\   \forall  h\in [0, L(Y(l))].
\end{gather}

Using induction with respect to $t$ it is not difficult to obtain from \eqref{m5} that for any $t>l$  and $h\in [0, L(Y(l))]$
$$     {\bf E}e^{hS(t)}\leq \sup_{0\le  t\le l}{\bf E}e^{hS(t)}.$$
Hence
\begin{gather}\label{m6}
\sup_{0\le  t\le l}{\bf E}e^{hS(t)}\le \sup_{t\ge 0}{\bf E}e^{hS(t)}\leq \sup_{0\le t\le l}{\bf E}e^{hS(t)}
,\ \forall  h\in [0, L(Y(l))].
\end{gather}

Then the results are proved from  (\ref{m6}), Theorem \ref{thm2}  and Corollary \ref{c2}.

\hfill\fbox

The following corollary is a special case of Theorem \ref{thm3}.
\begin{cor}\label{C-4}
Suppose  there exists a real  number $l>0,$  such that for  all $t\geq0$  claim sizes $Z(t+l)$ and  $Z(t)$  are identically distributed.  And  also  assume   that   (\ref{p19}) holds.
Then the  assertions  of   Theorem  \ref{thm3}  still hold.
\end{cor}

Model satisfying  conditions of Corollary \ref{C-3} or Theorem \ref{thm3}   will be called  quasi-periodic. Model  satisfying conditions  of Corollary \ref{C-4}  with
$$ \forall t\ge0, \  e^{-r(t+l)}\Lambda^{'}(t+l)= e^{-r(t)}\Lambda^{'}(t) \text{ and }
e^{-r(t+l)}p^{'}(t+l)= e^{-r(t)}p^{'}(t),$$
is  naturally  be called   periodic (or even purely periodic).

\begin{rem}
Let  $r(t)\equiv 0$   in  Corollary \ref{C-3}   or  Theorem  \ref{thm3} then  it  is  not  difficult  to  see  that the  similar   results  also
hold  for  periodic  risk  models  without  interest  rates.
\end{rem}

\begin{rem}
 From   Theorem \ref{thm3} we find several evident advantages of our results  compared    with the results in Asmussen and  Rolski \cite{AR94}. Firstly, in Theorem \ref{thm3} if we let  $r(t)=0$ and, let   (\ref{p18})  and (\ref{p19}) be equalities, then  we can obtain the results which were  presented in \cite{AR94}. Secondly, from  Remark \ref{r1}  and  (\ref{p18+})  it is easy  to see  that our estimate is  better. Thirdly, we do not  use the relatively difficult method  by a  change of measure which needs  technical details.  Our proof is based on a martingale approach  with less assumptions.
\end{rem}

\begin{exa}  \label{exa-0}
Suppose  the  period  $l=2$  and the  claim  sizes $Z(t)$  have   exponential  distribution    $f(x)=e^{-x},$  $ x>0.$
We  also  assume  that  $p^{'}(t)=4t$   and $\Lambda^{'}(t)=1,$  then they  satisfy  all  the conditions of  Theorem   \ref{thm3}  with $r(t)=0.$

 We  can  calculate  that  for  $h<1$
  $${\bf E}e^{hZ(t)}=\int_{0}^{\infty}e^{hx}e^{-x}dx=\frac{1}{1-h}. $$
 Then for any $t>0$
  \begin{gather*}
   {\bf E}e^{hS(t)}=\exp\Big\{-h\int_{0}^t 4xdx+\int_{0}^t \frac{h}{1-h}dx      \Big\}\\
 =\exp\Big\{ -2ht^2+th/(1-h)     \Big\}.
 \end{gather*}
 It  is  easy  to see  that
 $$L(S(2))=\sup\Big\{h\ge 0:   {\bf  E}e^{hS(2)}\le 1     \Big\}=3/4.$$
 Then  by  Theorem \ref{thm3}, for any $u>0$   we have
 $$    \psi(u)\le    \inf_{h\in [0, 3/4]} \big\{   e^{-hu}\sup_{0\le t\le 2}{\bf E}e^{hS(t)}\big\}\le \frac{3}{2}e^{-\frac{3}{4}u}.$$
\end{exa}

\subsection{  United risk model}
Insurance companies can reduce the probability of ruin by investing their assets in risk-free assets or    risky  assets. In  addition,  some insurance  companies  may transfer   part of the claims to the reinsurance company  to reduce  the  probability  of  ruin. More  details  about  the  investment  and  reinsurance can  be  found  in  Schmidli \cite{S07} and, Bai and Guo \cite{BG08}.

Here  we  consider  an  interesting  risk  model in  which   a number  of claim surplus processes are incorporated  such that the ith claim surplus process is added   at time $t^{(i)},$ which means that at time $t^{(i)}$ the insurance company will get  another premium income with new premium rate and may need to pay new claims.
Here we call these claim surplus processes  as a number of branches of the insurance company such that branch $i$  begins to work at time $t^{(i)}.$

Let $S^{(0)}(t),$ $S^{(1)}(t), S^{(2)}(t), \ldots$  be a sequence of   independent  claim  surplus  processes  such that
$$\forall t\ge 0,\      S^{(i)}(t)=\sum_{k=1}^{N^{(i)}(t)}Z^{(i)}(T^{(i)}_k)-p^{(i)}t,\  i=0, 1, \ldots,\text{ with }  S^{(i)}(0)=0,$$
where $N^{(i)}(t)$  is  the  Poisson  process  with  intensity  $\lambda^{(i)},$ $p^{(i)}$  is  the  non-random  premium   rate  and  $\{Z^{(i)}(T^{(i)}_k),  k=0,1, 2, \ldots\}$  are  i.i.d.  claim  sizes   which  are  independent  of  $N^{(i)}(t)$   where
$$T^{(i)}_k:=\min \{t\ge 0: N^{(i)}(t)\ge k\},\ k=0,1 ,\ldots.$$
From  (\ref{d9}) that for   each $  i=0, 1, 2, \ldots$ and any $h\ge 0$
\begin{gather}\label{u1}
{\bf E}e^{hS^{(i)}(t)}=e^{ta^{(i)}(h)}\text{ with } a^{(i)}(h)=\lambda^{(i)}({\bf E}e^{hZ^{(i)}}-1)-hp^{(i)}.
\end{gather}
And  let  $0=t^{(0)}<t^{(1)}< \ldots$  be a sequence of non-random times.   Define  a new   claim  surplus  process
\begin{gather}                           \label{u2}
S(t):=\sum_{i=0}^{\infty}S^{(i)}((t-t^{(i)})^+).
\end{gather}
We  may interpret    $S^{(i)}(t)$   as a process describing the \emph{ith} class of business of the insurance company
 which begins  to work  at time $t^{(i)}.$

Thus   we  consider  the  following  surplus  process
\begin{gather}\label{0a3}
   R(t)=u-S(t),\  \forall    t\ge 0,
\end{gather}
with  the  initial  surplus  $R(0)=u>0$  and  $S(t)$  defined  in  (\ref{u2}).
In  this  case, we  say  that   the  insurer's surplus  $R(t)$   defined  in  (\ref{0a3})  varies  according   to  the united  risk  model.

We  can  see  that   this  model can also be studied by  applying  model  A.

For any $t>0$  and $h\ge 0$ introduce notations
\begin{gather}\label{u3}
a(h,t)=\sum_{i=0}^{\infty}  (t-t^{(i)})^+a^{(i)}(h),\
a_{k}(h)=\sum_{i=1}^k (t^{(k)}-t^{(i)})a^{(i)}(h), \ a(h)=\sup_{k\geq 0}a_k(h).
\end{gather}
And also introduce ontations
\begin{gather}\label{u4}
\overline{L}:=\sup\big\{h\geq0 : a(h) \leq 0     \big\}, \ L^{(i)}: =\sup\big\{h\geq0 : a^{(i)}(h)\leq 0      \big\}.
\end{gather}

\begin{cor}\label{C-6}
Under the above  statements, suppose  $S(t)$ defined  in  (\ref{u2}) is  a  separable  process  with  independent  increments then  the Lundberg   inequality   in  (\ref{a2}) holds   with  $L:=\overline{L}$
where
\begin{gather}\label{u5}
L^{(0)}\geq \overline{L}\geq \inf_{i\geq 0}L^{(i)}.
\end{gather}
\end{cor}

{\bf  Proof.} First, it  is  clear  that (\ref{u5})  holds.
By (\ref{u1}), (\ref{u2})  and the property  of independency,  we get  that
\begin{eqnarray*}
{\bf E}e^{hS(t)}&=&\prod_{i=0}^{\infty} {\bf E}e^{hS^{(i)}((t-t^{(i)})^+)}\nonumber\\
&=&\prod_{i=0}^{\infty} e^{(t-t^{(i)})^+a^{(i)}(h)}\nonumber\\
&=&e^{\sum_{i=0}^{\infty}  (t-t^{(i)})^+a^{(i)}(h)   }\nonumber\\
&=&e^{a(h,t)}.
\end{eqnarray*}

By  (\ref{u3})  it is obvious  to obtain that
\begin{eqnarray*}
\sup_{t\geq 0}a(h,t)=\sup_{k\geq 0}a_k(h),
\end{eqnarray*}
since $a(h, \cdot)$   is a linear piecewise function.

The rest  of the proof  comes from   Theorem   \ref{thm-1}, we omit it here.

\hfill\fbox

Next  we   are  going  to present   an  example such  that     $\overline{L}=L^{(0)}.$
\begin{exa} \label{exa-1}
Suppose that $L^{(0)}>0$  and
\begin{gather*}
a^{(i)}(L^{(0)})>0,\ i=1, 3, 5, 7, \ldots,
\end{gather*}
\begin{gather}\label{u7}
a^{(i-1)}(L^{(0)})+a^{(i)}(L^{(0)}) \leq 0, \ i=1, 3, 5, 7, \ldots.
\end{gather}
It may be interpreted that claim surplus processes  with numbers $1, 3, 5, 7, \ldots$   have negative income.  But before opening "bad"  branch we open "good"  branch with  property  (\ref{u7}). In this case $\overline{L}=L^{(0)}.$
\end{exa}

\begin{rem}
If  $i=0$  in   the  united  risk  model then    it  degenerates  into  a  homogeneous compound Poisson  risk  model
with
$${\bf E}e^{hS(t)}=exp\big\{t[\lambda^{(0)}({\bf E}e^{hZ^{(0)}}-1)-hp^{(0)}] \big\}.$$
\end{rem}

\section{Comparisons with the existing results}

In this section, we are going to continue studying the result of non-homogeneous renewal risk model from Zhou, Sakhanenko ang Guo \cite{ZSG19} and comparing with the results in Andrulyt$\dot{e}$ et  al. \cite{AE15},    Kievinait$\dot{e}$ and $\check{S}$iaulys \cite{KD18}, and Kizinevi$\check{c}$ and  $\check{S}$iaulys  \cite{EJ18}.

Consider a class of risk processes $R(t)$ with the following properties:

\hspace*{0.5cm} (1) Process $R(t)$ may have positive jumps only at random or non-random times $T_1, T_2, \ldots$ such that
 \begin{gather*}
\forall k=1, 2, \ldots, \ T_{k+1}>T_{k} >T_0:=0 \text{ and } T_n\rightarrow\infty\ a.s.
\end{gather*}

\hspace*{0.5cm} (2) Process $R(t)$ is monotone on each interval $[T_{k-1}, T_k), \ k=1,2,\ldots,$ and
$R(0)=u>0.$

{\bf Model R.}
 Assume  the $k$-th  claim $Z_k$ occurs at  time $T_k$, i.e.,
\begin{gather*}                                                               
-Z_k:= R(T_k)-R(T_k-0)\le 0,\  \theta_k:=T_k-T_{k-1}>0,\  k=1,2,\ldots.
\end{gather*}
Suppose that   on each interval $[T_{k-1}, T_k)$  the premium  rate is $p_k,$ i.e.,
\begin{gather*}
\forall t\in[T_{k-1}, T_k),\  R(t)-R(T_{k-1})=p_k(t-T_{k-1}),\   k=1,2,\ldots.
\end{gather*}
Assume also  that random vectors
\begin{gather*}\label{aa7}
(p_k, Z_k, \theta_k),\   k=1,2,\ldots
\end{gather*}
are mutually independent.

Then conditions (1) and (2) hold and  the   random variables
\begin{gather}\label{aa3}
Y_k=R(T_{k-1})-R(T_k)=Z_k-p_k \theta_k=Z_k-X_k,\  k=1,2,\ldots
\end{gather}
are also mutually independent in  which   $X_k:=p_k \theta_k$  is the nonnegative accumulated  premium over a period of   time   $\theta_k.$ Model R is also called non-homogeneous renewal risk model which is different from the model A introduced in section  \ref{sec2.1}.
The   ruin  probability   of    model R  is   that   for any $u>0$
$$ \psi(u)={\bf P}(\sup_{k\ge 1}S_k >u),$$
where  $S_k=Y_1+Y_2+\ldots+Y_k.$    Theorem  1  of    \cite{ZSG19} states that for any $u>0$ and $h\ge 0$
\begin{gather}\label{z30}
  \psi(u)\le e^{-hu}\sup_{n\ge 1}{\bf E}e^{hS_n}.
\end{gather}


Then, we  are  going to  estimate  ${\bf E}e^{hS_n}$  in  (\ref{z30}) to obtain a  further explicit upper bound of ruin probability  of non-homogeneous renewal risk model and compare our estimation with the results in Andrulyt$\dot{e}$ et  al. \cite{AE15},    Kievinait$\dot{e}$ and $\check{S}$iaulys \cite{KD18}, and Kizinevi$\check{c}$ and  $\check{S}$iaulys  \cite{EJ18}.

The following simple assertion will be useful below.

\begin{lem}\label{L2}
Let random variables $Y_1, Y_2, \ldots$ have finite expectations. Then for any $h\ge 0$  and $m\ge 1$
\begin{gather}                                                                \label{z31}
 M_m(h):=\max_{1\leq k\leq m}{\bf E}e^{hS_k}\leq e^{h C_m}{\bf E}e^{hS_m},
\end{gather}
where
\begin{gather}                                                                \label{z32}
0\le C_m:=\max_{1\leq k\leq m}{\bf E}S_k-{\bf E}S_m<\infty.
\end{gather}
\end{lem}
{\bf Proof.} For any $1\le k\le m$  and $h\ge 0,$ by the Jensen's inequality in  probability space we have
\begin{gather}\label{A1}
{\bf E}e^{h S_n}={\bf E}e^{hS_k}\cdot{\bf E}e^{h(S_m-S_k)}\ge {\bf E}e^{h S_k}\cdot e^{h({\bf E}S_m-{\bf E}S_k)}
\geq
{\bf E}e^{hS_k}\cdot e^{-hC_m}.
\end{gather}

From (\ref{A1}) we immediately have
\begin{gather}\label{A2}
{\bf E}e^{hS_k}\le e^{hC_m}{\bf E} e^{hS_m}.
\end{gather}
By taking maximum on both sides of (\ref{A2}) over $1\le k\le m$ the result is obtained.

\hfill\fbox


The following inequality is evident and  useful for us: for any $h\ge 0$  and $1\le m\le n$
\begin{gather}\label{za4}
M(h):=\sup_{n\geq 1}{\bf E}e^{hS_n}\le M_m(h)\cdot\big(\sup_{n>m}{\bf E}e^{h(S_n-S_m)} \vee 1 \big).
\end{gather}
In fact, for any $h\ge 0$ and $1\le m\le n$
\begin{eqnarray*}
M(h):&=&\sup_{n\ge 1}{\bf E}e^{hS_n}\\
&=&(\sup_{n\le m}{\bf E}e^{hS_n})\vee (\sup_{n>m}{\bf E}e^{hS_n})\\
&=&M_m(h)\vee \{ \sup_{n>m}{\bf E}e^{h(S_n-S_m)}\cdot e^{hS_m}   \}\\
&=& M_m(h) \vee \{ {\bf E}e^{hS_m}\cdot \sup_{n>m}{\bf E}e^{h(S_n-S_m)}  \}\\
&\le& M_m(h)\cdot\{1\vee \sup_{n>m}{\bf E}e^{h(S_n-S_m)}   \},
\end{eqnarray*}
since ${\bf E}e^{hS_m}\le M_m(h).$

The estimations  of  $M(h)$  and $M_n(h)$ are the central  technical problems in  applications  of  (\ref{z30}). Many ideas are known  in this direction.  The  most famous ones are connected  with papers of Bennett \cite{BG62}
 and Hoeffding \cite{HW63}. In this  section, we are going to recall  and use the ideas from \cite{BG62}
 and  \cite{HW63}. First, by using  the important limits of mathematical analysis such that $\lim_{n\rightarrow \infty}(1+\frac{1}{n})^n=e$ we have
 \begin{gather*}                                                                    
 {\bf E}e^{hS_n}=\prod_{k=1}^{n}{\bf E}e^{hY_k}\leq \Big(\frac{1}{n}\sum_{k=1}^{n}{\bf E}e^{hY_k}\Big)^{n}\leq \exp\big\{ \sum_{k=1}^{n}({\bf E}e^{hY_k}-1)   \big\}.
 \end{gather*}
Next, introduce non-negative function $g(x)=e^x-1-x\geq 0$ and note that
$${\bf E}e^{hY_k}=1+h{\bf E}Y_k+{\bf E}g(hY_k).$$
 From  \cite{BG62} and \cite{HW63} we see that function $g(x)$   has many  useful   properties. We  are going to use them here.
As we all know for each $k=1,2,\ldots$
\begin{gather*}                                                                    
 Y_k=Z_k-X_k \text{ with } X_k, Z_k\geq 0,
\end{gather*}
then
\begin{gather*}                                                                  \label{z5}
{\bf E}g(h Y_k)={\bf E}[g(h Y_k):Y_k\geq0]+{\bf E}[g(h Y_k):Y_k<0]\\
= {\bf E}[g(hY^{+}_k)]+ {\bf E}[g(-hY^{-}_{k})]
\leq {\bf E}g(h Z_k)+{\bf E}g(-h X_k).
\end{gather*}
It is easy to see that for any $h,C\geq 0$
\begin{gather*}                                                                  
 {\bf E}g(-h X_k)\leq \frac{h^2}{2}{\bf E}\big[ X^{2}_k:0<X_k \leq C   \big]+ h {\bf E}\big[X_k: X_k>C  \big].
\end{gather*}
But the most remarkable for our purpose  is the following property:
if $ H>0$   by using the error term of  Taylor expansion of $e^{x}$  at $x=0$  such that
$\delta(x)=e^x-1-x=x^2 \int^1_0 (1-t) e^{tx}dt$   and the definition of $g(x)$  with $x=hZ_k$ we have
\begin{eqnarray}                              \label{z7}
\forall  h\in(0,H],\
h^{-2}{\bf E}g(h Z_k)&=&{\bf E}\int_0^1Z_k^2(1-t)e^{thZ_k}dt\nonumber\\
&\leq&
{\bf E}\int_0^1Z_k^2(1-t)e^{tHZ_k}dt=H^{-2}{\bf E}g(H Z_k).
\end{eqnarray}

For $H,C>0$ introduce notations
\begin{gather}                                                                    \label{z8}
A_n(C)=\sum_{k=1}^{n}{\bf E}\big[X_k: X_k>C   \big],
\ B_n(H, C)=\frac{1}{H^2}\sum_{k=1}^{n}{\bf E}g(HZ_k)+\frac{1}{2}\sum_{k=1}^{n}{\bf E}\big[ X_{k}^{2}: X_k\leq C  \big].
\end{gather}
From the above statements we can see that
\begin{gather}\label{z9}
{\bf E}e^{hS_n}\leq \exp\big\{  \sum_{k=1}^{n}\big[h {\bf E} Y_k+ {\bf E}g(h Y_k)\big]  \big\}
\leq  \exp\big\{h {\bf E} S_n +h A_n(C)+h^2 B_n (h, C)   \big\}.
\end{gather}

Then  from  \eqref{z9} we have  the following assertion.
\begin{cor}\label{cor8}
Suppose that for some integer  $m>0$ and  real numbers  $h,C>0,$
\begin{gather}\label{z10}
A_n(C)+ hB_n(h, C)\leq -{\bf E}S_n,\ \forall n\ge m.
\end{gather}
Then
\begin{gather}\label{z11}
M(h)\leq e^{h C_m},
\end{gather}
where  $C_m$ is  defined in \eqref{z32}.
In particular, for any $u>0$
\begin{gather}                                                                    \label{z12}
\psi(u)\leq   e^{h(C_m-u)}.
\end{gather}
\end{cor}
{\bf Proof.} From (\ref{z9})  and (\ref{z10}), it is easy to see that
\begin{gather}\label{za1}
{\bf E}e^{hS_n}\le \exp \{ h{\bf E}S_n+h(-{\bf E}S_n) \}=1,\  \forall  n\ge m.
\end{gather}
Then by inequalities   of  (\ref{z31}), (\ref{za4}) and (\ref{za1})  we obtain that
\begin{eqnarray*}
M(h)&\le& M_m(h)(\sup_{n>m}{\bf E}e^{h(S_n-S_m)}\vee 1)\\
&\le& e^{hC_m}{\bf E}e^{hS_m}(\sup_{n>m}{\bf E}e^{h(S_n-S_m)}\vee 1)\\
&=&e^{hC_m}({\bf E}e^{hS_m}\vee \sup_{n>m}{\bf E}e^{hS_n})\\
&=&e^{hC_m}\cdot \sup_{n\ge m}{\bf E}e^{hS_n}\\
&\le& e^{hC_m}.
\end{eqnarray*}

(\ref{z12}) is directly  obtained from (\ref{z30})  and (\ref{z11}).

\hfill\fbox

If  (\ref{z10})  does not  hold, the following fact will be useful.
\begin{cor}\label{cor9}
Suppose that  for some integer $m>0$
\begin{gather}                                                                    \label{z13}
 A_n(C)\leq -c^{*}{\bf E}S_n   \text{ with } c^*\in (0,1),  \  \forall n\ge m,
\end{gather}
and
\begin{gather}                                                                    \label{z14}
 B_n(H, C)\leq C^*(-{\bf E}S_n) \text{ with } 0<H,C^*<\infty,\  \forall n\geq m.
\end{gather}
Then all assertions of
 Corollary  \ref{cor8} hold  with
\begin{gather}                                                                    \label{z15}
0<h\le\min\{H, (1-c^*)/C^{*}\}.
\end{gather}
\end{cor}

{\bf Proof.}
It follows from (\ref{z7}) and (\ref{z8}) that $B_n(h, C)\leq B_n(H, C)$ for $ h\in(0,H]$.
From this fact and  (\ref{z14}), for any $ h\in(0,H]$ we have
\begin{gather}                              \label{z15+}
hB_n(h, C)\leq hB_n(H, C)\leq hC^*(-{\bf E}S_n).
\end{gather}
Thus,  inequality (\ref{z15+}), with  $h$ in (\ref{z15}), yields that $hB_n(h, C)\leq (1-c^*)(-{\bf E}S_n)$
and, hence, (\ref{z10}) is satisfied in view of (\ref{z13}).

\hfill\fbox

\begin{cor}\label{cor10}
Suppose that
\begin{gather}\label{z16}
\limsup_{n\rightarrow \infty}{\bf E}S_n<0,\
\lim_{C\rightarrow \infty}\limsup_{n\rightarrow \infty}\frac{A_n(C)}{-{\bf E}S_n}<1,
\end{gather}
and that for some $H>0$  and any $C>0$
\begin{gather}\label{z17}
 \limsup_{n\rightarrow \infty}\frac{B_n(H, C)}{-{\bf E}S_n}<\infty.
\end{gather}
Then there exists $h>0$  for which (\ref{z11})  and (\ref{z12}) hold.
\end{cor}
{\bf Proof.} It follows from  (\ref{z16}) that there exists an integer $m>0$ such that
$\sup_{n\ge m}{\bf E}S_n<0$ and that (\ref{z13}) holds with some  $c^*\in(0,1)$.
After that we have  from (\ref{z17}) that (\ref{z14}) is evidently true with some $C^*<\infty$.
So, Corollary \ref{cor10} follows from Corollary \ref{cor9}.

\hfill\fbox

\begin{rem}
It is clear that
\begin{gather*}
{\bf E}g(HZ_k)\leq {\bf E}e^{HZ_k}-1\leq {\bf E}e^{HZ_k},
\end{gather*}
since $H>0$ and $Z_k$ is nonnegative random variable.
And  for all $C>0$
\begin{gather}                                                                    \label{z19}
{\bf E}[X^{2}_k: X_k\leq C]\leq C^2 {\bf P}[X_k>0]\leq C^2.
\end{gather}
It is not difficult to verify that from  these facts  and Corollaries \ref{cor9} and \ref{cor10}   all  assertions in \cite{AE15, KD18, EJ18}  can be obtained.

Moreover, in this way we may obtain  better constants  than  in  \cite{AE15, KD18, EJ18}, since in  these three  papers    the   authors used that
\begin{gather*}\                                                                    \label{z20}
{\bf E}g(HZ_k/2)\leq \frac{1}{2}(\frac{H}{2})^2{\bf E}\big[Z^{2}_k e^{HZ_k/2}   \big]
\leq 2e^{-2} {\bf E}e^{H Z_k}
\end{gather*}
instead of  (\ref{z7}).
\end{rem}

Note that we may always put $Z_k=Y_k^+$ and $X_k=Y_k^-$.
In this case (\ref{z19}) is not trivial with ${\bf P}[X_k>0]={\bf P}[Y_k<0]$.

\begin{exa}\label{exa3}
Let independent random variables  $Y_1, Y_2, \ldots$  have normal distributions such that
\begin{gather*}                                                                                 
B_n:={\bf Var} S_n\to \infty,
\\                                                                                   
\forall  k=1, 2, \ldots,\ {\bf E}Y_k:=-(27/64)(B_k+B_{k-1}){\bf Var} Y_k.
\end{gather*}
(We may, for simplicity,  take $B_k:=k$.)
In this case ${\bf E}S_n=-(27/64)B_n^2$ and
$${\bf E}e^{hS_n}=e^{h {\bf E}S_n+h^2\var S_n/2}= e^{-(27/64)hB_n^2+h^2B_n/2}=e^{f(B_n,h)}.$$
It is easy to see that
\begin{gather*}\label{exa3.1}
\sup_{n\geq 1} {\bf E}e^{hS_n}=\sup_{n\geq 1}e^{f(B_n,h)}
\leq\sup_{x\geq 0}e^{f(x,h)}
=e^{f(16h/27,h)}=e^{4h^3/27}.
\end{gather*}
Hence
$$e^{-hu}\sup_{n\geq1}{\bf E}e^{hS_n}\leq e^{-hu+4h^3/27}=e^{g(h,u)}.$$
Then for any $u>0$
$$\psi(u)\le\inf_{h\geq0}e^{-hu}\sup_{n\geq1}{\bf E}e^{hS_n}
\le\inf_{h\geq0}e^{g(h,u)}=e^{g((3/2)\sqrt u,u)}=e^{-u^{3/2}}.
$$

Thus, it is possible in non-homogeneous case that we use $h=h(u)\rightarrow\infty$
as $u\to\infty$ to obtain better bound than $\psi(u)=e^{-O(u)}$. Such situation is impossible in the case of homogeneity.
\end{exa}

\bigskip

{\noindent {\bf\large Acknowledgments.}  The authors are thankful to the referees for careful reading
of the paper and helpful comments and suggestions.

\flushleft Qianqian Zhou \hfill Alexander Sakhanenko\\
 Nankai University  \hfill Sobolev Institute of Mathematics\\
School of Mathematical Sciences  \hfill 4 Acad. Koptyug avenue\\
Tianjin 300071, China  \hfill  Novosibirsk 630090, Russia.\\
  E-mail: qianqzhou@yeah.net \hfill Nankai University \\
  \hfill School of Mathematical Sciences\\
\hfill Tianjin 300071, China\\
\hfill  Email: aisakh@mail.ru

\flushleft Junyi Guo\\
 Nankai University \\
School of Mathematical Sciences\\
Tianjin 300071, China\\
  E-mail: jyguo@nankai.edu.cn


\begin{thebibliography}{1234}

\bibitem{AE15} I. M. Andrulyt$\dot{e}$, E.  Bernackait$\dot{e}$, D.  Kievinait$\dot{e}$ and J. $\check{S}$iaulys,
A Lundberg-type inequality for an inhomogeneous renewal risk model, Modern Stochastics: Theory and Applications 2(2)(2015), pp. 173-184.
\bibitem{AA10} S. Asmussen and H. Albrecher,   Ruin Probabilities, World Scientific Publishing 2010.
\bibitem{AR94} S.  Asmussen and T. Rolski,   Risk theory in a periodic environment: the Cramer-Lundberg approximation and Lundberg's inequality,    Mathematics of Operations Research 19(2)(1994), pp.410-433.
\bibitem{BG08}  L. Bai and  J. Guo.   Optimal proportional reinsurance and investment with multiple risky assets and no-shorting constraint. Insurance Mathematics and Economics, 42(3)(2008) pp.968-975.
 \bibitem{BG62} G. Bennett,  Probability inequalities for the sum of independent random variables,  Journal of the American Statistical Association,  57 (1962), pp. 33-45.
 \bibitem{BJ17} E. Bernackait$\dot{e}$  and J. $\check{S}$iaulys,  The finite-time ruin probability for an inhomogeneous renewal risk model, Journal of Industrial and Management Optimization,  13(1)(2017) pp.207-222.
 \bibitem{BBE10} K. Bla$\check{z}$evi$\check{c}$ius, E.  Bieliauskien$\dot{e}$  and J.  $\check{S}$iaulys, Finite-time ruin probability in the inhomogeneous claim case. Lithuanian Mathematical Journal,  50(3)(2010) pp.260-270.
 \bibitem{CCGLM13}  A. Casta$\tilde{n}$er, M. M.  Claramunt,  Gathy M., et al, Ruin problems for a discrete time risk model with non-homogeneous conditions. Scandinavian Actuarial Journal, 2013(2)(2013) pp.83-102.
\bibitem{C30} H. Cram$\acute{e}$r,   On the mathematical theory of risk, Skandia Jubilee. Vol. 4, Stockholm 1930.
\bibitem{C55} H. Cram$\acute{e}$r,    Collective  risk  theory, Jubilee Volume, Skandia Insurance Company 1955.
\bibitem{DC05} D. C. M. Dickson,   Insurance Risk and Ruin, Cambridge: Cambridge University Press 2005.
\bibitem{GH79} H.  Gerber,   An  introduction  to  Mathematical  Risk Theorey, S.S. Huebner  Foundation  Monographs, University  of  Pennsylvania 1979.
\bibitem{G91} J. Grandell,   Aspects of Risk Theory,   Springer, New York 1991.
\bibitem{HW63} W. Hoeffding,  Probability inequalities for sums of bounded random variables, J. Am. Stat.
Assoc.  58(301) (1963) pp. 13-30.
\bibitem{IK00} Z. G. Ignatov and V. K.  Kaishev, Two-sided bounds for the finite-time probability  of  ruin. Scandinavian  Actuarial  Journal,  2000(1)(2000) pp.46-62.
\bibitem{KD18} D. Kievinait$\dot{e}$ and J.  $\check{S}$iaulys,  Exponential bounds for the tail probability of the supremum of an inhomogeneous random walk, Modern Stochastics: Theory and Applications 5 (2) (2018), pp.129-143.
\bibitem{KJF93} J. F. C. Kingman,   Poisson  processes.  Clarendon  Press, Oxford 1993.
\bibitem{EJ18} E. Kizinevi$\check{c}$ and J. $\check{S}$iaulys,   The Exponential Estimate of the Ultimate Ruin Probability for the Non-Homogeneous Renewal Risk Model, Risks 6(1)(2018), 20.
\bibitem{LP06} C. Lef$\grave{e}$vre and P. Picard,  A nonhomogeneous risk model for insurance, Computers and Mathematics with Applications, 51(2)(2006) pp.325-334.
\bibitem{LF03} F. Lundberg,  Approximerad  framst$\ddot{a}$ llning  av  sannolikhetsfunktionen.  ${\AA}$  terf$\ddot{o}$rs$\ddot{a}$kring av kolletivrisker. Acad. Afhaddling. Almqvist. och  Wiksell, Uppsala 1903.
\bibitem{PJ97}  J. Paulsen,    Present value of some insurance portfolios. Scandinavian Actuarial Journal 1997(1) (1997) pp.11-37.
\bibitem{RSST98} T. Rolski, H. Schmidli, V.  Schmidt and J. L. Teugels,  Stochastic Processes for Insurance and Finance,
John Wiley, Chichester 1998.
\bibitem{S07} H.  Schmidli,  Stochastic Control in Insurance, Springer-Verlag, London 2007.
\bibitem{TT14} A. Tuncel and  F. Tank, Computational results on the compound binomial risk model with nonhomogeneous claim occurrences. Journal of Computational and Applied Mathematics,  263(2014) pp.69-77.
\bibitem{VR15} R. Vernic,   On a conjecture related to the ruin probability for nonhomogeneous insurance claims, Annals of the Ovidius University of Constanta-Mathematics Series,  23(3)(2015) pp.209-220.
\bibitem{ZSG19}  Q.Q. Zhou,  A. Sakhanenko and J. Y.  Guo,  Lundberg-type inequalities  for non-homogeneous risk models, arXive: 2004.11190v1, 2020.
\end{thebibliography}
\end{document}